\documentclass[10pt]{amsart}
\usepackage{amssymb}
\usepackage{amsfonts}
\usepackage{amscd}
\usepackage[dvips]{graphicx}
\usepackage[all,cmtip]{xy}
\usepackage{hyperref}
\usepackage{yhmath}
\usepackage{bbm}
\usepackage{ulem}
\usepackage[all]{xy}
\usepackage{mathdots}
\usepackage{leftidx}
\usepackage{mathrsfs}
\usepackage{amsmath,amssymb, amsthm}
\usepackage{color}
\usepackage{tikz}
\usepackage{mathdots}
\usepackage{picture}
\usepackage{enumerate}
\usepackage{makecell}
\usepackage{extarrows}
\usepackage{graphicx}
\usepackage{caption}
\usepackage{subcaption}
\usepackage{float}
\numberwithin{equation}{section}
\newtheorem{Them}{Theorem}[section]
\newtheorem{Lem}[Them]{Lemma}
\newtheorem{Def}[Them]{Definition}
\newtheorem{Cor}[Them]{Corollary}
\newtheorem{Prop}[Them]{Proposition}
\newtheorem{Ex}[Them]{Example}
\newtheorem{Rem}[Them]{Remark}

\newtheorem{Prob}[Them]{Problem}

\newcommand{\add}{{\mathsf{add}}}
\newcommand{\Add}{{\mathsf{Add}}}
\newcommand{\Prod}{{\mathsf{Prod}}}

\newcommand{\End}{{\mathsf{End}}}

\newcommand{\Id}{{\mathsf{Id}}}

\newcommand{\DTr}{{\mathsf{DTr}}}
\newcommand{\Tr}{{\mathsf{Tr}}}
\newcommand{\D}{{\mathsf{D}}}

\newcommand{\m}{\mathsf{mod}}

\newcommand{\stmod}{ \mathsf{\underline{mod}}}
\newcommand{\Stmod}{ \mathsf{\overline{mod}}}
\newcommand{\Hom}{{\mathsf {Hom}}}
\newcommand{\Coker}{{\mathsf {Coker}}}
\newcommand{\Ker}{{\mathsf {Ker}}}
\newcommand{\Img}{{\mathsf {Im}}}
\newcommand{\StHom}{\mathsf{\underline{Hom}}}

\title[Simple-minded systems and coherent rings]{Simple-minded systems and coherent rings}
\author{Zhen Zhang}

\address{Zhen Zhang
	\newline College of Education for the future
	\newline Beijing Normal University 
	\newline  Zhuhai 519087
	\newline P.R.China}
\email{zhangzhen@bnu.edu.cn}

\date{version of \today}

\newcommand{\sdp}{\times\kern-.2em\vrule height1.1ex depth-.05ex}

 \normalbaselines
\begin{document}
	\renewcommand{\thefootnote}{\alph{footnote}}
	\renewcommand{\thefootnote}{\alph{footnote}}
	\setcounter{footnote}{-1} 
	\renewcommand{\thefootnote}{\alph{footnote}}
	\setcounter{footnote}{-1}

\begin{abstract}
Let $A$	be a  finite dimensional algebra over an algebraically closed field. We present a relationship between simple-minded systems and coherent rings.
\end{abstract}

\maketitle
\section{Introduction}
Simple-minded systems (sms for short), introduced by Koenig and Liu \cite{KL}, is a family  of objects which satisfies orthogonality and generating condition over the stable module category of any artin algebra.  Dugas \cite{Dugas} defined simple-minded systems in any Hom-finite Krull-Schmidt triangulated category. The authors \cite{KL,CKL,CLZ,GLYZ,Z} studied properties of sms's, including the finite cardinality of an sms,  invariance of an sms under stable equivalences, characterization of an sms over a  representation-finite self-injective algebra and the location of the members of an sms on a quasi-tube over a self-injective algebra and so on.  
In this paper, we generalize sms's to a right (resp. left) triangulated category (refer to \cite{BM1,KV}) which is a natural generalization of an abelian category, the stable module category of any artin algebra, or a triangulated category.  Note that Coelho Sim\~{o}es \cite{C} introduced $d$-simple-minded systems over $(-d)$-Calabi-Yau triangulated categories for any positive integer $d$ and there  is a recent rise of interests in studying $d$-simple-minded systems (see, for example, \cite{CP,IJ}).

Weakly simple-minded systems (wsms for short) was  introduced by Koenig-Liu \cite{KL} to compere with simple-minded systems, and they showed that a family of objects $\mathcal{S}$ in $A$-$\stmod$ is an sms if and only if  $\mathcal{S}$ is a wsms over a representation-finite self-injective algebra $A$.   We \cite{Z} present an example that there is a wsms which is not an sms over the stable module category of a $2$-domestic Brauer graph algebra.  In this paper, we shall study further relationship between simple-minded systems and weakly simple-minded systems over a right (resp. left) triangulated category. 
	
	
Let $\mathcal{T}$ be a  triangulated category and  $\mathcal{X}$ a subset of an sms. Dugas \cite{Dugas} showed that $(^{\perp}\mathcal{X}, \mathcal{F}(\mathcal{X}))$ and $(\mathcal{F}(\mathcal{X}),\mathcal{X}^{\perp})$ are torsion pairs and $\mathcal{F}(\mathcal{X})$ is functorially finite in $\mathcal{T}$. We know that for a wsms $\mathcal{S}$ in $\mathcal{T}$, if $\mathcal{F}(\mathcal{S})$ is covariantly finite or contravariantly finite in $\mathcal{T}$, then  $\mathcal{S}$ is an sms. It indicates that an sms  has a close relationship with functorially finiteness for an orthogonal system.  We now ask one question as follows.
\begin{Prob}
For an orthogonal system  $\mathcal{X}$ over a right {\rm(}resp. left{\rm)} triangulated category $\mathcal{X}$,  when  the extension  closure $\mathcal{F}(\mathcal{X})$ is covariantly finite or contravariantly finite?
\end{Prob}


It is easy to see that if  $\mathcal{T}$ is the stable module (projective) category of any representation-finite algebra, then $\mathcal{F}(\mathcal{X})$ is always functorially finite.  However, it is a hard problem in a general setting, such as the stable module (projective) categories of  representation-infinite  algebras.
Let $R$ be a ring and $\mathcal{M}$ a class of modules  in $R$-module. 
It was shown by Enochs \cite{E} that the category of all finitely generated projective modules is covariantly finite in $R$-module if and only if $R$ is right coherent. And more generally, it was proved by Crawley-Boevey \cite{CB} that covariantly finiteness of $\mathcal{M}$ in  $R$-module is characterized by closure under products of the category $\varinjlim \mathcal{M}$. It means that the existence of  covariantly finiteness of $\mathcal{M}$  is related to coherence. 
In this paper, by the coherence and approximation theory, we  explore  relationship between simple-minded systems and weakly simple-minded systems over a right (resp. left) triangulated category. 


This paper is organized as follows. In Section 2, we recall right triangulated category, approximation theory, simple-minded system and related definitions and conclusions. In Section 3, we prove main results Theorem \ref{properties-subset-of-sms}, Theorem \ref{fin-dim-sms} and Theorem \ref{weak-sms}.

	
	\section{Preliminaries}	
Throughout this paper, we assume that all additive categories are Krull-Schmidt, that is, any object is isomorphic to a finite direct sum of objects whose endomorphism rings are local. We assume that all rings are assumed to be associative with identity and all algebras  are considered to be finite dimensional algebras  over algebraically closed field $k$. For a ring $R$, we denote by $R$-$\m$ the category of finite presented left $R$-modules. The stable module category of $R$-$\m$ is denoted by $R$-$\stmod$, it has the same class of objects with $R$-$\m$, and for two objects $X,Y$ in $R$-$\stmod$, the abelian group   $\StHom_R(X,Y)$ from $X$ to $Y$ is the quotient  $\StHom_R(X,Y)/\mathcal{P}(X,Y)$, where $\mathcal{P}(X,Y)$ is the subgroup of $\StHom_R(X,Y)$ consisting of all $R$-module homomorphisms  which factor through a projective  $R$-module.  
	
We denote by $\Omega$ (resp. $\Sigma$) the syzygy functor (resp. cosyzygy functor) which assigns to any object $M$ of $R$-$\stmod$ the kernel of its projective cover $P_{R}(M)\twoheadrightarrow M$ (cokernel of its injective enveloping $M\hookrightarrow I_{R}(M)$) in $R$-$\m$. We  denote by $\add(M)$ (resp. $\Add(M)$) all direct summands of finite (resp. arbitrary) direct sums of modules from $M$, and $\Prod(M)$ all direct summands of arbitrary direct sums of modules from $M$.  Note that, for a class of objects $\mathcal{M}=\{M_{i}\mid i\in I\}$ in an additive category $\mathcal{C}$, we simply denote $\add(\bigoplus_{i\in I}M_{i})$ by $\add(\mathcal{M})$. In particular, $\mathcal{P}_{R}$ (resp. $\mathcal{I}_{R}$) is denoted by the full subcategories of finite generated  projective (resp. injective) $R$-modules $\add(_{R}R)$ (resp. \!$\add(\D(R_{R}))$). 
In this paper, when dealing with endomorphism rings we will use the convention of writing endomorphism opposite scalars, that is, if $f$ and $g$ are elements of an endomorphism ring, then $fg=g\circ f.$


\subsection{ Right triangulated categories}
Let $\mathcal{C}$ be an additive category equipped with an additive endofunctor $\Sigma:\mathcal{C}\rightarrow\mathcal{C}$. Consider the category $\mathcal{RT}(\mathcal{C},\Sigma)$: its objects are of the form
\[ \xymatrix{  X\ar[r]^-{} & Y\ar[r]^-{} &  Z\ar[r]^-{} & \Sigma(X), }\] 
and the morphisms are the triple $(f,g,h)$ of morpisms
indicated by the following diagram:
\[\xymatrix{
X \ar[d]_{f} \ar[r]^{\alpha} & Y \ar[d]_{g} \ar[r]^{\beta} & Z \ar[d]_{h} \ar[r]^{\gamma}&	\Sigma(X) \ar[d]^{\Sigma(f)}  \\
X' \ar[r]^{\alpha'} & Y' \ar[r]^{\beta'} & Z'  \ar[r]^{\gamma'}  &\Sigma(X').} \]
The composition of morphisms of $\mathcal{RT}(\mathcal{C},\Sigma)$ is induced in the canonical way by the corresponding composition of the morphisms in $\mathcal{C}$.

A right triangulation of the pair $(\mathcal{C}, \Sigma)$ is a full subcategory $\nabla$ of $\mathcal{RT}(\mathcal{C},\Sigma)$ which satisfies all the axioms of a triangulated category, except that $\Sigma$ is not  necessarily an equivalence. Then the triple  $(\mathcal{C},\Sigma,\nabla)$, is called a  right triangulated category.  We state the definition of a right triangulated category as follows.

\begin{Def}\label{left-right-tiangulated} {\rm(\cite[Chapter II, Section 1, Definition 1.1]{BR})} 
A full subcategory $\nabla$ of	$\mathcal{RT}(\mathcal{C},\Sigma)$ is said to be a {\bf right triangulation} of  $(\mathcal{C},\Sigma)$ if it is closed under isomorphisms and satisfies the following four axioms:
	\begin{enumerate}[$(1)$]
		\item For any object $X$ of $\mathcal{C}$, the right triangle  $\xymatrix{ 0 \ar[r]^-{} &   X\ar[r]^-{1_{X}} &X\ar[r]^-{} &  0}$ belongs to $\nabla$.
		\item For any morphism $h:X\rightarrow Y$, there is a right triangle in $\nabla$ of the form \[\xymatrix{   X\ar[r]^-{} & Y\ar[r]^-{} &  Z\ar[r]^-{}&\Sigma(Z).}\]

		\item For any two right  triangles \[\xymatrix{   X\ar[r]^-{f} & Y\ar[r]^-{g} &  Z\ar[r]^-{h}&\Sigma(X)  }\ \ and\ \  \xymatrix{   X'\ar[r]^-{f'} & Y'\ar[r]^-{g'} &  Z'\ar[r]^-{h'}&\Sigma(X') }\] 
		in $\nabla$, and for any morphism
		 $\alpha:X\rightarrow X'$ and  $\beta:Y\rightarrow Y'$ of $\mathcal{C}$ with $\beta\circ f=f'\circ\alpha$, there is a morphism  $\gamma:Z\rightarrow Z'$ such that $(\alpha,\beta,\gamma)$ is a morphism from the first triangle to the second, see the following diagram:
		 \[\xymatrix{
		 X \ar[d]_{\alpha} \ar[r]^{f} & Y \ar[d]_{\beta} \ar[r]^{g} & Z \ar@{-->}[d]_{\gamma} \ar[r]^{h} & 	\Sigma(X)\ar[d]^{\Sigma(\alpha)} \\
		 	X' \ar[r]^{f'} & Y' \ar[r]^{g'} & Z'\ar[r]^{h'}  & \Sigma(X').} \]

		\item ${\rm(}${\bf Octahedral axiom}${\rm)}$ For any two right  triangles \[\xymatrix{   X\ar[r]^-{f} & Y\ar[r]^-{g} &  Z\ar[r]^-{h}&\Sigma(Z)  }\ \ and\ \  \xymatrix{   Y\ar[r]^-{f'} & M\ar[r]^-{g'} &  N\ar[r]^-{h'}&\Sigma(Y) }\] 
		in $\nabla$, there is a third right triangle $\xymatrix{   X\ar[r]^-{f'\circ f} & M\ar[r]^-{u} &  U\ar[r]^-{v}&\Sigma(X)  }$ in $\nabla$, and two morphism  $s:Z\rightarrow U$ and  $t:U\rightarrow N$ of $\mathcal{C}$  such that the diagram below is fully commutative and 
		\[\xymatrix{   Z\ar[r]^-{s} & U\ar[r]^-{t} &  N\ar[r]^-{(\Sigma g)\circ h'}&\Sigma(Z)  }\] is a right triangle in $\nabla$.
		 \[\xymatrix{
			X \ar@{=}[d]_{\Id_{X}} \ar[r]^{f} & Y \ar[d]_{f'} \ar[r]^{g} & Z \ar@{-->}[d]_{s} \ar[r]^{h} & 	\Sigma(X)\ar@{=}[d]^{\Sigma(\Id_{X})} \\
			X \ar[r]^{f'\circ f} & M\ar[d]_{g'} \ar[r]^{u} & U\ar@{-->}[d]_{t}\ar[r]^{v}  & \Sigma(X) \ar[d]^{\Sigma f}\\
				& N\ar[d]_{h'}\ar@{=}[r]^{\Id_{N}} & N\ar[d]^{(\Sigma g)\circ h'}\ar[r]^{h'}  & \Sigma(Y) \\
				 & \Sigma Y \ar[r]^{\Sigma g} & \Sigma Z &  \\
			} \]
\end{enumerate}
The triple $(\mathcal{C},\Sigma,\nabla)$  is called a {\bf right triangulated category}. The functor $\Sigma$ is called  {\bf suspension functor} and a diagram in $\Sigma$ is called  a {\bf right triangle}. 
\end{Def} 

{\bf Left triangulated category} is defined dually, denoted by the triple  $(\mathcal{C},\Omega,\Delta)$, $ \Omega$ is called {\bf loop functor} and the diagrams {\rm(}of the form $\xymatrix{\Omega(Z)\ar[r] & X\ar[r]^-{} & Y\ar[r]^-{} &  Z}${\rm)} in $\Delta$ are called {\bf left triangles}. Note that the octahedral axiom holds for left (resp. right) triangulated categories. Please refer \cite{BM1, BM2, BR, KV} for basic properties of a  left (resp. right) triangulated category.

\begin{Rem}
\begin{enumerate}[$(1)$]
\item Triangulated categories are both left and right  triangulated categories, here $\nabla=\Delta$ and $\Omega=\Sigma^{-1}$.

\item Any additive category with kernels {\rm(}resp. cokernels{\rm)} is a left {\rm(}resp. right{\rm)}  triangulated category {\rm(}in particular any abelian category{\rm)}, where $\Omega=\Sigma=0$, and $\Delta$ is the class of left exact sequences and $\nabla$ the class of right exact sequences. 

\end{enumerate}
\end{Rem}
 Let $\mathcal{C}$ be an additive category and $\mathcal{X}$ a full subcategory of $\mathcal{C}$.  A morphism $f: A\rightarrow B$ is said to be an {\bf$\mathcal{X}$-epic} (resp. {\bf$\mathcal{X}$-monic}), if the morphism $\mathcal{C}(X,f): \mathcal{C}(X,A)\rightarrow \mathcal{C}(X,B)$ (resp. $\mathcal{C}(f,X):\mathcal{C}(B,X)\rightarrow\mathcal{C}(A,X)$) is surjective for each $X\in\mathcal{X}$.

We say that a morphism $f\colon M\rightarrow X$ is a {\bf left $\mathcal{X}$-approximation} of $M$ if $X\in\mathcal{X}$, and the abelian group homomorphism $\mathcal{C}(f,X')\colon \mathcal{C}(X,X')\rightarrow \mathcal{C}(M,X')$ is surjective for each $X'\in\mathcal{X}$. Dually, a morphism $g\colon Y\rightarrow N$ is a {\bf right $\mathcal{X}$-approximation} of $N$ if $Y\in\mathcal{X}$, and the abelian group homomorphism $\mathcal{C}(Y',g)\colon \mathcal{C}(Y',Y)\rightarrow \mathcal{C}(Y',N)$ is surjective for each $Y'\in\mathcal{X}$. Let $\mathcal{S}$ be a full subcategory of $\mathcal{C}$ containing $\mathcal{X}$. $\mathcal{S}$ is said to be {\bf covariantly finite} (resp. {\bf contravariantly finite}) in $\mathcal{X}$ if every $S\in\mathcal{S}$ has a left (resp. right) $\add(\mathcal{X})$-approximation. We say $\mathcal{S}$  is  {\bf functorially finite} in  $\mathcal{X}$ provided that  $\mathcal{S}$ is both covariantly finite and contravariantly finite in $\mathcal{X}$. 
 

It is known that if additive category $\mathcal{C}$ has coproduct (resp. product), and $\mathcal{X}\subseteq\mathcal{C}$ is skeletally small, then $\Add(\mathcal{X})$(resp. $\Prod(\mathcal{X})$) is contravariantly finite (resp. covariantly finite) in $\mathcal{C}$.

According to Beligiannis and Marmaridis \cite{BM2}, we have the following general conclusion. Let $\mathcal{X}$ (resp. $\mathcal{Y}$) be a covariantly finite (resp. contravariantly finite) subcategory of  $\mathcal{C}$ and assume that $\mathcal{C}$ is abelian, or more generally that  any $\mathcal{X}$-monic has a cokernel  (resp. any $\mathcal{Y}$-epic has a kernel) in $\mathcal{C}$, then $\mathcal{C}/\mathcal{X}$ admits a right triangulated structure $(\mathcal{C}/\mathcal{X}, \Sigma_{\mathcal{X}}, \nabla_{\mathcal{X}})$ (resp. a left triangulated structure $(\mathcal{C}/\mathcal{Y}, \Omega_{\mathcal{Y}}, \Delta_{\mathcal{Y}})$), where $\Sigma_{\mathcal{X}}$ (resp. $\Omega_{\mathcal{Y}}$) is the suspension functor (resp. the loop functor) and $\nabla_{\mathcal{X}}$ (resp. $\Delta_{\mathcal{Y}}$) is the right  triangulation (left triangulation). We list examples as follows. 

\begin{Ex}$($\cite[chapter II, Section 1, p. 24$\sim$26]{BR}$)$ \label{ex-pretri}
Let $R$ be an artin ring. It is known that both $\mathcal{P}_{R}$ and $\mathcal{I}_{R}$ are functorially finite in  $R$-$\m$. Then 
\begin{enumerate}[$(1)$]
\item The stable module {\rm(}projective{\rm)} category $R$-$\stmod:=R$-$\m/\mathcal{P}_{R}$ is a left triangulated category with loop functor the syzygy functor $\Omega$, where $\mathcal{P}_{R}$ is the full subcategory of finitely generated projective left $R$-modules.  Note that $R$-$\stmod$ is also a right triangulated category with  suspension functor $\Sigma=\Tr\Omega\Tr$, where $\Tr$ is the transpose duality functor. 
\item The stable module {\rm(}injective{\rm)} category $R$-$\Stmod:= R$-$\m/\mathcal{I}_{R}$ is a right triangulated category with suspension functor the cosyzygy functor $\Sigma$, where $\mathcal{I}_{R}$ is the full subcategory of finitely generated injective left $R$-modules. Note that $R$-$\Stmod$ is also a left triangulated category with loop functor $\Omega=\DTr\Omega\Tr\D$, where $\Tr$ is the transpose duality functor and $\D$ is the usual duality  for algebra $R$.  
\end{enumerate}
\end{Ex}

\subsection{Simple-minded systems}
Based on  Dugas' definition of extension closure \cite[Section 2]{Dugas} on a triangulated category, we have the following routine generalization on a right (resp. left) triangulated category. Let $(\mathcal{T},\Sigma,\nabla)$ be a right triangulated category.   For any families $\mathcal{S}_{1}, \mathcal{S}_{2}$ of objects in $\mathcal{T}$, we define a family of objects
\[\mathcal{S}_{1}\star\mathcal{S}_{2}:=\{ X\in \mathcal{T}\mid \mbox{ There is a right triangle }S_{1} \longrightarrow  X \longrightarrow S_{2} \longrightarrow \Sigma S_{1}\in\nabla, where \ S_{1}\in \mathcal{S}_{1}, S_{2}\in \mathcal{S}_{2}\}. \]
	
 Using the octahedral axiom,  the operator $\star$ is associative, that is,  $(\mathcal{S}_{1}\star\mathcal{S}_{2})\star\mathcal{S}_{3}=\mathcal{S}_{1}\star(\mathcal{S}_{2}\star\mathcal{S}_{3})$ for $\mathcal{S}_{1}, \mathcal{S}_{2}$ and $\mathcal{S}_{3}\subseteq \mathcal{T}$.  We denote $(\mathcal{S})_{0}=\{0\}$, and for $n\in\mathbb{Z}^{+}$, we inductively define $(\mathcal{S})_{n}=(\mathcal{S})_{n-1}\star(\mathcal{S}\cup\{0\})$. We have $(\mathcal{S})_{n}\star(\mathcal{S})_{m}=(\mathcal{S})_{n+m}$ for any non-negative integers $m$ and $n$. Similarly, one can define $ _{n}(\mathcal{S})$, and we have $(\mathcal{S})_{n}$=$_{n}(\mathcal{S})$.  We say that $\mathcal{S}$
is {\it extension-closed}, if $\mathcal{S}\star\mathcal{S}\subseteq \mathcal{S}$. One denotes the {\bf extension closure} of a family $\mathcal{S}$ of objects in $\mathcal{T}$ as $$\mathcal{F}(\mathcal{S}):=\bigcup_{n\geq0}(\mathcal{S})_{n},$$ which is the smallest extension closed full subcategory of  $\mathcal{T}$  containing $\mathcal{S}$. Dually, one may define extension closure over a left triangulated category. 
	
For any family $\mathcal{S}$ of objects in $\mathcal{T}$, we set 
$$\mathcal{S^{\perp}}:=\{Y\in \mathcal{T} \mid\mathcal{T}(X,Y)=0, \forall X\in \mathcal{S}\},$$
$$\mathcal{^{\perp}S}:=\{Y\in \mathcal{T}\mid\mathcal{T}(Y,X)=0, \forall X\in \mathcal{S}\}.$$ 
We know that both $\mathcal{S^{\perp}}$ and $\mathcal{^{\perp}S}$ are extension closed subcategories of $\mathcal{T}$  as well as closed under direct summands. We shall denote $\mathcal{S^{\perp}}\cap\mathcal{^{\perp}S}$ by $\mathcal{^{\perp}S^{\perp}}$, which is called a {\it stable bi-perpendicular category}. Note that $(\mathcal{^{\perp}S})^{\perp}$ (resp. $^{\perp}(\mathcal{S}^{\perp})$) has a different meaning with stable bi-perpendicular category $\mathcal{^{\perp}S^{\perp}}$, $\mathcal{S}$ is contained in $(\mathcal{^{\perp}S})^{\perp}$ (resp. $^{\perp}(\mathcal{S}^{\perp})$), but not contained in $\mathcal{^{\perp}S^{\perp}}$.

\begin{Def}\label{brick-orthogonal-system}
Let $\mathcal{T}$ be an additive $k$-category.  An object $M$ in $\mathcal{T}$ is a {\bf stable brick} if $\mathcal{T}(M,M)\cong k$.   Moreover, a family $\mathcal{S}$ of stable bricks in $\mathcal{T}$  is an {\bf orthogonal system} if $\mathcal{T}(M,N)=0$ for all distinct  $M, N$ in $\mathcal{S}$.
\end{Def}

\begin{Rem}\label{F-S-properties}
It is routine to check that, if $\mathcal{S}$  is an orthogonal system, then $\mathcal{F}(\mathcal{S})$ is  closed under direct summands. The reader may refer to \cite[Lemma 2.7]{Dugas} for more details.
\end{Rem}

\begin{Def}\label{torsion-pair} 
Let $(\mathcal{T},\Sigma,\nabla)$ be a right triangulated category and let $\mathcal{X}, \mathcal{Y}$ be two additive subcatgories of $\mathcal{T}$ which are closed under direct summands and isomorphisms. The pair $(\mathcal{X},\mathcal{Y})$ is called a {\bf torsion pair}, if 
\begin{enumerate}[$(1)$]
\item $\mathcal{T}(\mathcal{X},\mathcal{Y})=0$.
\item For any $C\in\mathcal{T}$, there a  right triangle		
\[\xymatrix{X_{C}\ar[r]^-{} & C\ar[r]^-{} &  Y^{C}\ar[r]^-{}& \Sigma(X_{C})}\]
with $X_{C}\in\mathcal{X}$ and $Y^{C}\in\mathcal{Y}$.
\end{enumerate}
\end{Def}

\begin{Rem}\label{homogeneous-tube-tir}
\begin{enumerate}[$(1)$]
\item The torsion pair in Definition \ref{torsion-pair} is slightly different from the torsion pair in \cite[Chapter II, Section 3, Definition 3.1]{BR}, the condition $\Sigma(\mathcal{X})\subseteq\mathcal{X}$ is not included in our definition. 
\item We know that, for a subcategory $\mathcal{X}$ of a right triangulated category $\mathcal{T}$, both $(\mathcal{X}^{\perp},{^{\perp}(\mathcal{X}^{\perp})})$ and $((^{\perp}\mathcal{X})^{\perp}, {^{\perp}\mathcal{X}})$ 
are torsion pairs.
\end{enumerate}
\end{Rem}

\begin{Def}\label{definition-sms-right-tir} {\rm(\cite[Definition 2.4, 2.5]{Dugas})} 
Let $(\mathcal{T},\Sigma,\nabla)$ be a right triangulated category.  A family of objects $\mathcal{S}$ in $\mathcal{T}$ is a {\bf simple-minded system} {\rm(sms for short)} if the following two conditions are satisfied$\colon$
\begin{enumerate}[$(1)$]
\item {\rm(Orthogonality)} $\mathcal{S}$ is an orthogonal system in $\mathcal{T}$. 
\item {\rm(Generating condition)} Extension closure $\mathcal{F}(\mathcal{S})$ of $\mathcal{S}$ is equal to $\mathcal{T}$.
\end{enumerate}
\end{Def}
	

\begin{Def}\label{definition-wsms-right-tir} {\rm(\cite[Definition 5.3]{KL})} 
Let $(\mathcal{T},\Sigma,\nabla)$ be a right triangulated category.  A family of objects $\mathcal{S}$ in $\mathcal{T}$ is a {\bf weakly simple-minded system} {\rm(wsms for short)} if the following two conditions are satisfied$\colon$
\begin{enumerate}[$(1)$]
\item {\rm(Orthogonality)} $\mathcal{S}$ is an orthogonal system in $\mathcal{T}$. 
\item {\rm(Weakly generating condition)} For any non-zero object $X\in\mathcal{T}$, there is an object $S$ in $\mathcal{S}$ such that $\mathcal{T}(X,S)\ncong 0.$
\end{enumerate}
\end{Def}

\begin{Rem}\label{wsms-sms}
\begin{enumerate}[$(1)$]
\item It follows from Koenig-Liu \cite{KL}  that simple-minded systems coincide with weakly simple-minded systems over the stable module category of any representation-finite self-injective algebra. 
\item We \cite{Z} present an example that there is a wsms which is not an sms over the stable module category of a $2$-domestic Brauer graph algebra. 
\end{enumerate}
\end{Rem}

We state two necessary conditions of sms's as follows.
\begin{Them}$($\cite[Theorem 3.3]{Dugas}$)$\label{subset-of-sms}
Let $\mathcal{T}$ be a Hom-finite Krull-Schmidt right {\rm(}or left{\rm)} triangulated category. Suppose $\mathcal{X}\subseteq\mathcal{S}$ for a simple-minded system $\mathcal{S}$ in $\mathcal{T}$. Then $(^{\perp}\mathcal{X},\mathcal{F}(\mathcal{X}))$ and $(\mathcal{F}(\mathcal{X}),\mathcal{X}^{\perp})$ are torsion pairs in $\mathcal{T}$.  In particular, $ \mathcal{F}(\mathcal{X})$ is a functorially finite subcategory of $\mathcal{T}$.
\end{Them}	
\begin{proof}
It is a direct generalization of \cite[Theorem 3.3]{Dugas}.
\end{proof}
\begin{Them}$($\cite[Theorem 1.3]{CLZ}$)$\label{simple-module-and-n-tube}
Let $A$ be a self-injective algebra and $\mathcal{C}$ a quasi-tube of rank $n$. Then the number of elements in an sms of $A$ lying in $\mathcal{C}$ is strictly less than $n$. In particular, none of the indecomposable modules in an sms  lie in the homogeneous tubes of the AR-quiver.
\end{Them}	


\subsection{Coherent  rings and direct limit}

Let $R$ be a ring. $R$ is said to be {\bf left (resp. right) coherent} if every finitely generated left (resp. right) ideal of $R$ is finitely presented. Recall that a left $R$-module $X$ is called {\bf finitely presented} if there is an exact sequence 
\[\xymatrix{ R^{n} \rightarrow R^{m} \rightarrow X\rightarrow 0}\]
 of $R$-modules, where $m,n$ are positive integers.  In general, let $M$, $X$ be two left $R$-modules, then $X$ is called {\bf finitely $M$-presented} if $X$ admitting an exact sequence \[\xymatrix{ M^{n} \rightarrow M^{m} \rightarrow X\rightarrow 0}\] with $m,n$ positive integers.  Note that finitely $M$-copresented modules can be defined dually.
 Recall that a left $R$-module $F$ is {\bf flat} if $-\otimes_{R}F$ is an exact functor, that is,  whenever 
\[\xymatrix{0\ar[r]^-{} & A\ar[r]^-{g} & B\ar[r]^-{h}& C\ar[r]^{}&0}\]
is an exact sequence of right $R$-modules, then 
\[\xymatrix{0\ar[r]^-{} & A\otimes_{R}F\ar[r]^-{g\otimes\Id_{B}} & B\otimes_{R}F\ar[r]^-{h\otimes\Id_{C}}& C\otimes_{R}F\ar[r]^{}&0}\]
is an exact sequence of abelian groups.

\begin{Them}$($\cite[Theorem 2.1]{Chase}$)$\label{flat-coherent}
For any ring  $R$ the following statements are equivalent.
\begin{enumerate}[$(1)$]
    \item  The direct product of any family of flat right $R$-modules is flat.
	\item The direct product of any family of copies of $R$ is flat as a right $R$-module.
	\item Any finitely generated submodule of a free left $R$-module is finitely presented.
	\item Any finitely generated left ideal in $R$ is finitely presented.
\end{enumerate}
\end{Them}

\begin{Rem}\label{coherent-rings-flat}
Let $R$ be a ring and $M$ a left $R$-module. 
\begin{enumerate}[$(1)$]
\item A  ring $R$ is  left coherent if and only if the direct product of any family of flat right $R$-modules is flat $($cf. \cite{Chase}$)$.	
\item There are left coherent rings which are not right coherent {\rm(}\cite[Example 4.46(e)]{Lam}{\rm)}.
\end{enumerate}
\end{Rem}


Let $R$ be a ring,  $M$ a left $R$-module and $\mathcal{Q}=(Q_{i}\mid i\in I)$  a sequence of left $R$-modules. Let $\varphi_{M,\mathcal{Q}}$ denote the abelian group homomorphism
\begin{equation}\label{prod-iso}
 \varphi_{M,\mathcal{Q}}\colon M\otimes_{R} \prod_{i\in I}Q_{i}\longrightarrow\prod_{i\in I}M\otimes_{R} Q_{i}.
 \end{equation}
\begin{equation*}
\ \ \ \ \ \ \ \ \  m\otimes(q_{i})_{i\in I}\mapsto (m\otimes q_{i})_{i\in I}
\end{equation*}
\begin{Lem}$($\cite[Chapter 3, Section 3.1, Lemma 3.8]{GT}$)$\label{epi-iso}
Let $R$ be a ring and $M$  a left $R$-module.
\begin{enumerate}[$(1)$]
\item	$M$ is finitely generated, if and only if  $\varphi_{M,\mathcal{Q}}$ is an epimorphism for all sequences $\mathcal{Q}$ of left $R$-modules, if and only if $\varphi_{M,\mathcal{Q}}$ is an epimorphism for the sequence $\mathcal{Q}=(Q_{i}\mid i\in I)$ such that $I=M$ and $Q_{i}=R$ for all $i\in I$.
\item $M$ is finitely presented, if and only if  $\varphi_{M,\mathcal{Q}}$ is an isomorphism for all sequences $\mathcal{Q}$ of left $R$-modules. 
\end{enumerate}
\end{Lem} 	

\begin{Cor}\label{coherent-epi-iso}
 Let $R$ be a left coherent ring and let $M$ be  finitely generated as a right $R$-module. If $\mathcal{Q}=(Q_{i}\mid i\in I)$ is  a sequence of flat left $R$-modules, then $\varphi_{M,\mathcal{Q}}$ in {\rm(}\ref{prod-iso}{\rm)} is an isomorphism.
\end{Cor}
\begin{proof}
Since $M$ is finitely generated, there is a short exact sequence
\begin{equation}\label{fin-gen}
\xymatrix{0\ar[r]^-{} & K\ar[r]^-{g} & F\ar[r]^-{h}& M\ar[r]^{}&0,}
\end{equation}
where $F$  is a finitely generated free left $R$-module and $K$ is a finitely generated left $R$-module. By Lemma \ref{epi-iso}, $\varphi_{K,\mathcal{Q}}$ is  epimorphic.
Since $R$ is  coherent and $Q_{i}$ is flat for each $i\in I$, the direct product $\prod_{i\in I}Q_{i}$ is flat as a left $R$-module. Tensoring with sequence (\ref{fin-gen}), we have 
 \[\xymatrix{
0\ar[r]&	K\otimes_{R}\prod_{i\in I}Q_{i} \ar[d]_{\varphi_{K,\mathcal{Q}}} \ar[r]^{g\otimes  \Id_{\prod_{i\in I}Q_{i}}}& F\otimes_{R}\prod_{i\in I}Q_{i} \ar[d]_{\varphi_{F,\mathcal{Q}}} \ar[r]^{h\otimes  \Id_{\prod_{i\in I}Q_{i}}} & M\otimes_{R}\prod_{i\in I}Q_{i}\ar[d]_{\varphi_{M,\mathcal{Q}}} \ar[r] & 	0 \\
0\ar[r]&\prod_{i\in I}K\otimes_{R}Q_{i} \ar[r]^{({g\otimes \Id_{Q_{i}}})_{i\in I}} & \prod_{i\in I}F\otimes_{R}Q_{i}\ar[r]^{({h\otimes \Id_{Q_{i}}})_{i\in I}}  & \prod_{i\in I}M\otimes_{R}Q_{i}\ar[r] & 0.} \]
	
Since $K$ and $M$ are finitely generated, by Lemma \ref{epi-iso}, both $\varphi_{K,\mathcal{Q}}$ and $\varphi_{M,\mathcal{Q}}$ are  epimorphic. Since $F$ is a finitely generated free module, then $\varphi_{F,\mathcal{Q}}$ is an isomorphism. By Snake lemma, there is an exact sequence as follows.
 
\[\xymatrix{0\ar[r]^-{} &\Ker(\varphi_{K,\mathcal{Q}}) \ar[r]^-{} &\Ker(\varphi_{F,\mathcal{Q}})\ar[r]^-{} 
& \Ker(\varphi_{M,\mathcal{Q}}) \ar[d]^-{} \\
& & &\Coker(\varphi_{K,\mathcal{Q}}) \ar[r]^-{}& \Coker(\varphi_{F,\mathcal{Q}}) \ar[r]^{}& \Coker(\varphi_{M,\mathcal{Q}}) \ar[r]^{}&0.}\]
 
 Since $\Ker(\varphi_{F,\mathcal{Q}})\cong \Coker(\varphi_{F,\mathcal{Q}})\cong 0$,
 $\Ker(\varphi_{M,\mathcal{Q}})\cong \Coker(\varphi_{K,\mathcal{Q}})\cong0$, thus $\varphi_{M,\mathcal{Q}}$ is an isomorphism.	
\end{proof}

\begin{Rem}
In general, $\varphi_{M,\mathcal{Q}}$ is not a monomorphism, and the condition is not  easily captured. $M$ is called a Mittag-Leffler module provided that $\varphi_{M,\mathcal{Q}}$ in {\rm(}\ref{prod-iso}{\rm)} is monomorphism for all sequences $\mathcal{Q}$ of left $R$-modules. Refer to \cite[Chapter 3]{GT} for related conclusions.
\end{Rem}

\begin{Def}\label{direct-limit} 
Let  $\mathcal{T}$ be a category, $I$  a partially ordered set,  and let $\{M_{i},\varphi_{j}^{i}\}$ be a direct system in  $\mathcal{T}$ over $I$. The {\bf direct limit} is an object $\varinjlim M_{i}$ and insertion morphism $(\alpha_{i}\colon M_{i}\rightarrow\varinjlim M_{i})_{i\in I}$ such that 
\begin{enumerate}[$(1)$]
\item $\alpha_{j}\circ\varphi_{j}^{i}=\alpha_{i}$ whenever $i\preceq j$. 
\item Let $X\in\mathcal{T}$, and let there be given morphisms $f_{i}\colon M_{i}\rightarrow X$ satisfying  $f_{j}\circ\varphi_{j}^{i}=f_{i}$ for all $i\preceq j$. There exist a unique morphism $\theta\colon\varinjlim M_{i}\rightarrow X$ making the diagram commute:
	\[\xymatrix{\varinjlim M_{i}\ar@{-->}[rr]^{\theta} &&X\\
& M_{i} \ar@{->}[ul]_{\alpha_{i}} \ar@{->}[ur]^-{f_{i}} \ar[d]^{\varphi_{j}^{i}}&  \\
& M_{j}  \ar@/^-3ex/[uur]_{f_{j}} \ar@/^3ex/[uul]^{\alpha_{j}} &.}\]
\end{enumerate}
\end{Def}

Let $R$ be a ring and $Y$ a left $R$-module. A submodule $X$ of  $Y$ is a {\bf pure submodule} ($X\subseteq_{\ast}Y$ for short), if for each finitely presented module $P$, the functor $\Hom_{R}(P,-)$ preserves the exactness of the short exact sequence 	
\begin{equation}\label{pure-exact}
\xymatrix{0\ar[r]^-{} & X\ar[r]^-{g} & Y\ar[r]^-{h}& X/Y\ar[r]^{}&0.}
\end{equation}
The embedding is then called a {\bf pure embedding}, and the sequence $(\ref{pure-exact})$ is called a  {\bf pure-exact sequence}. An epimorphism $\pi:Y\rightarrow Z$ is called a {\bf pure epimorphism} provided that $\Ker\pi\subseteq_{\ast}Y.$
A $R$-module $Z$ is said to be {\bf pure-injective} if it is injective over pure embeddings, that is, if whenever $f : X\rightarrow Y$ is a pure embedding and $g : X\rightarrow Z$ is any morphism, then there is a morphism $h : Y\rightarrow Z$ with $g=h\circ f$, see the following diagram.
\[\xymatrix{
	X \ar[dr]_-{g} \ar[r]^-{f} & Y \ar@{-->}[d]^-{h} \\
	& Z.}\]

\begin{Lem}$($\cite[Chapter 2, Lemma 2.13]{GT}$)$\label{direct-lim-equiv}
Let $R$	be a  ring and $\mathcal{C}$ a class of finitely presented left $R$-modules closed under finite direct sums. Then the following are equivalent for a left $R$-module $M$.
\begin{enumerate}[$(1)$]
\item  $M\in\varinjlim\mathcal{C}$.
\item There is a pure epimorphism $f:\bigoplus _{i\in I}C_{i}\rightarrow M$ for a sequence $(C_{i}\mid i\in I)$ of modules in $\mathcal{C}$.
\item Every homomorphism $h:F\rightarrow M$, where $F$ is finitely presented, has a factorization through one module in $\mathcal{C}$.
\end{enumerate}
Moreover, $\varinjlim\mathcal{C}$ is closed under direct limit, pure submodules and pure epimorphic images, and the finitely presented modules in $\varinjlim\mathcal{C}$ are exactly the direct summands of modules in $\mathcal{C}$. 
\end{Lem}


The following result characterizes the class of  the direct limit of $\add( M)$ for an arbitrary module $M$.

\begin{Prop}$($\cite[Theorem 3.3]{PPT}$)$\label{direct-lim-of-add} 
Let $R$ be a  ring,  $M$  a left $R$-module and  $S=\End(_{R}M)$. Then $\varinjlim\add( M)$  coincides with the class of all modules of the form $M\otimes_{S}F$, where $F$ is a flat left $S$-module. 
\end{Prop} 

We recall some concepts as follows.
We say that a morphism $f:U\rightarrow V$ of $\add(M)$ is a {\bf pseudokernel} of morphism $g:V\rightarrow W$ when $g\circ f=0$ and if there is a morphism $h:U'\rightarrow V$ such that  $g\circ h=0$, then $h$ factors through $f$.  A left $R$-module $M$ is called an {\bf endoflat} module, if $M_{S}$ is a flat module, where $S$ is the endomorphism ring of $_{R}M$. 
Recall that a left (resp. right) $R$-module $M$ is {\bf finendo} if $M$ is finitely generated as a module over its endomorphism ring $\End({_{R}}M)$ (resp. $\End(M_{R})$). Dually, a left (resp. right) $R$-module $M$ is {\bf cofinendo} if $M$ is finitely cogenerated as a module over its endomorphism ring $\End({_{R}}M)$ (resp. $\End(M_{R})$). 
Recall that {\bf weak global dimension} $wD(S)$ is the supremum of projective dimensions of the finitely presented left $S$-modules.


\begin{Them}$($\cite[Proposition 1.1 and Theorem 1.8]{GMG}$)$\label{coherent-ring-add-kernel-1} 
Let $S$ be a  ring. The following conditions are equivalent.
\begin{enumerate}[$(1)$]
	\item Every morphism  of $\add(M)$ has a pseudokernel. 
	\item  Every finitely $M$-cogenerated module has a right  $\add(M)$-approximation.
	\item $S$  is left coherent.
\end{enumerate}
\end{Them} 

\begin{Them}$($\cite[Proposition 1.4 and Theorem 1.8]{GMG}$)$\label{coherent-ring-add-kernel-2} 
Let $S$ be a  ring. The following conditions are equivalent.
\begin{enumerate}[$(1)$]
		\item Every morphism of $\add(M)$ has a kernel. 
		\item $S$  is left coherent and weak global dimension $wD(S)$ of $R$ is at most $2$.
\end{enumerate}
\end{Them} 

\begin{Them}$($\cite[Theorem 1.6]{GMG}$)$\label{coherent-ring-add} 
	Let $R$ be a  ring,  $M$  a left $R$-module and  $S=\End(_{R}M)$. Then
	$S$  is left coherent, $M$ is endoflat and $wD(S)\leq2$ if and only if every  finitely $M$-copresented module belongs to $\add(M)$.
\end{Them} 



	
\section{Sms's and coherent rings}	

Beligiannis and Reiten \cite[Chapter V, Section 1, Lemma 1.1]{BR} present the following useful observation. For the convenience of the reader, we state a proof as follows.

\begin{Lem}\label{contra-variant}
Let $\mathcal{T}$ be an abelian category and let $\mathcal{X}\subseteq\mathcal{T}$ be a subcategory of $\mathcal{T}$. Assume that $\omega\subseteq\mathcal{X}$ and $\omega$ is a contravariantly finite {\rm(}resp. covariantly  finite{\rm)} in $\mathcal{T}$. Then $\mathcal{X}$ is contravariantly finite {\rm(}resp. covariantly  finite{\rm)} in $\mathcal{T}$ if and only if  $\mathcal{X}/\omega$ is contravariantly finite {\rm(}resp. covariantly  finite{\rm)} in $\mathcal{T}/\omega$.

\end{Lem} 
\begin{proof} We only prove the case that $\omega$ is a covariantly finite subcategory, since the other case can be proved dually.
Assume that  $\mathcal{X}$ is  covariantly  finite in $\mathcal{T}$ and let $M\in\mathcal{T}$. If $f_{M}: M\rightarrow X_{M}$ is a left $\mathcal{X}$-approximation of $M$ in $\mathcal{T}$, then it is easy to know that $\underline{f_{M}}: \underline{M}\rightarrow \underline{X_{M}}$ is a left $\mathcal{X}/\omega$-approximation of $M$ in $\mathcal{T}/\omega$. Thus,  if $\mathcal{X}$ is covariantly  finite in $\mathcal{T}$,  then $\mathcal{X}/\omega$ is  covariantly  finite in $\mathcal{T}/\omega$.
	
Conversely, let $\underline{\alpha}: \underline{M}\rightarrow \underline{X}$ be  a left $\mathcal{X}/\omega$-approximation of $M$ in $\mathcal{T}/\omega$. We can choose an object $X_{M}$ and a morphism $f_{M}: M\rightarrow X_{M}$ in $\mathcal{T}$ satisfying $\underline{X_{M}}=\underline{X}$ and $\underline{f_{M}}=\underline{\alpha}$. Since $\omega$ is  covariantly  finite in $\mathcal{T}$, without loss of generality, we may assume that $f_{M}$ is $\omega$-monic  (recall that a morphism $f:A\rightarrow B$ is said to be  $\mathcal{X}$-monic provided that the morphism $\mathcal{T}(f,\mathcal{X}):\mathcal{T}(B,\mathcal{X})\rightarrow\mathcal{T}(A,\mathcal{X})$ is surjective).
	
Let  $\beta: M\rightarrow X'$ be a morphism with $X'\in\mathcal{X}$. Since $\mathcal{X}/\omega$ is  covariantly  finite in $\mathcal{T}/\omega$, there is a morphism $\gamma: X_{M}\rightarrow X'$ such that
\begin{equation}\label{eq-1}
\underline{\gamma}\circ\underline{\alpha}=\underline{\gamma}\circ\underline{f_{M}}=\underline{\beta}.
\end{equation}	
See the following diagram:
\[\xymatrix{
	\underline{M} \ar[dr]_-{\underline{\beta}} \ar[r]^-{\underline{\alpha}=\underline{f_{M}}} & \underline{X_{M} } \ar@{-->}[d]^-{\underline{\gamma}}\\
&   \underline{X'}.}    \]       
	
Thus $\beta-\gamma\circ f_{M}$ factors through a left $\omega$-approximation $\varphi_{M}:M\rightarrow W$, that is, there is a morphism $\phi:W\rightarrow X'$ such that 
\begin{equation}\label{eq-2}
	\beta-\gamma\circ f_{M}=\phi\circ\varphi_{M}.
\end{equation}
Since $f_{M}$ is $\omega$-monic and $\omega\subseteq\mathcal{X}$, there is a morphism $\mu_{M}: X_{M}\rightarrow W$ such that \begin{equation}\label{eq-3}
	\mu_{M}\circ f_{M}=\varphi_{M}.
\end{equation}
See the following diagram:
\[\xymatrix{
&    M\ar[d]_-{\varphi_{M}} \ar[r]^-{\beta-\gamma\circ f_{M}} \ar[dl]_-{f_{M}}& X'  &  \\
	X_{M}\ar@{-->}[r]^-{\mu_{M}} & W \ar@{-->}[ur]^-{\phi}& & 
}   \]   
        
\medskip

Then we have $\beta\stackrel{(\ref{eq-2})}=\gamma\circ f_{M}+\phi\circ\varphi_{M}\stackrel{(\ref{eq-3})}=\gamma\circ f_{M}+\phi\circ\mu_{M}\circ f_{M}=(\gamma +\phi\circ\mu_{M})\circ f_{M}$, thus $\beta$ factors through $f_{M}$. Hence $f_{M}$ is a left $\mathcal{X}$-approximation of $\mathcal{T}$,  and  then $\mathcal{X}$ is  covariantly  finite in $\mathcal{T}$.
\end{proof}

\begin{Prop}$($\cite[Chapter 8, Lemma 8.35]{GT}$)$\label{direct-prod-definable}
Let $R$	be a left coherent ring and $\mathcal{C}$ a class of finitely presented modules  such that $\mathcal{C}=\add(\mathcal{C})$. Then the following are equivalent. 
\begin{enumerate}[$(1)$]
	\item  $\mathcal{C}$  is covariantly finite in $R$-$\m$.
	\item $\varinjlim\mathcal{C}$ is closed under direct products.
\end{enumerate}
\end{Prop}
\begin{proof}
$(1)$	implies $(2)$: Let $(C_{i}\mid i\in I)$ be a sequence of modules in $\mathcal{C}$, and take $M=\prod_{i\in I}C_{i}$. We show that $M$ belongs to $\varinjlim\mathcal{C}$. Let $f: F\rightarrow M$ be a homomorphism with $F$ finitely presented. Let $g: F\rightarrow C$ be a left $\mathcal{C}$-approximation  of $F$. Then for each $i\in I$, there is $h_{i}:C\rightarrow C_{i}$ such that $h_{i}\circ g=\pi_{i}\circ f$, where $\pi_{i}: M\rightarrow C_{i}$ is the canonical projection. Thus there is a morphism $h:C\rightarrow M$ such that $h_{i}=\pi_{i}\circ h$ for all $i\in I$, then $f=h\circ g$ (see the following diagram).  By Lemma \ref{direct-lim-equiv}, $M\in\varinjlim\mathcal{C}$.
\[\xymatrix{
	F \ar[d]_-{g} \ar[r]^-{f} & M \ar[dr]_-{\pi_{i}}   \\
	C \ar@{-->}[rr]_-{h_{i}}\ar@{-->}[ur]_{h} & & C_{i} 
}\]
	
\medskip
	
$(2)$	implies $(1)$:  Let $M$ be in $R$-$\m$. By Lemma \ref{direct-lim-equiv}, $\varinjlim\mathcal{C}$ is closed under pure submodules. Let $\kappa=|M|+|R|+ \aleph_{0}$ and $\mathcal{R}$ a representative set of all $\kappa$-generated modules in $\varinjlim\mathcal{C}$. Consider the canonical map $f:M\rightarrow L,$ where $L=\prod_{N\in\mathcal{R}}N^{\Hom_{R}(M,N)}\in\varinjlim\mathcal{C}$. Then any morphism $g:M\rightarrow N$ with $N\in\mathcal{R}$ factors through $f$.  
	
We claim that any morphism $f': M\rightarrow L'$  with $L'\in\varinjlim\mathcal{C}$ factors through $f$. 
Firstly, it is easy to see that $|\Img f|\leqslant \kappa$, so by \cite[Lemma 2.25 (a)]{GT}, there is a pure submodule $N'$ of $L'$ such that $\Img f\subseteq N'$ and $|N'|\leqslant\kappa$. Since $\varinjlim\mathcal{C}$ is closed under pure submodules, $N'\in\varinjlim\mathcal{C}$, and then  $N'\cong N$  for some $N\in\mathcal{R}$. Hence there is a morphism $\alpha:L\rightarrow L' $ such that $f'=\alpha\circ f$,  that is, $f'$ factors through $f$, see the following diagram:
\[\xymatrix{
M \ar[d]_-{f'}\ar[dr]^-{f'} \ar[r]^-{f} & L \ar@{-->}[d]^{\alpha} \\
	\ \ \ \ \ \ \ Im f' \hookrightarrow N' \ar[r]^-{\subseteq_{\ast}\ \ } & L'.
}\]

\medskip	
	
Since $M$ is finitely presented and $L\in\varinjlim\mathcal{C}$, by Lemma \ref{direct-lim-equiv} (2), there are $C\in\mathcal{C}$, $g: M\rightarrow C$ and $h: C\rightarrow L$, such that $f=h\circ g.$ Let $g': M\rightarrow C',$ where $C'\in\mathcal{C}$.  By the above claim, there is a morphism $h': L\rightarrow C'$ such that $g'=h'\circ f$, therefore $g'=h'\circ h\circ g$ (see the following diagram). Thus $g$ is a left $\mathcal{C}$-approximation  of $M$ in $R$-$\m$. 
\[\xymatrix{
	& C\ar[dr]^-{h}  \\
	M \ar[ur]^-{g} \ar[rr]^{f} \ar[dr]_-{g'} && L\ar@{-->}[dl]^{h'} \\
	& C' &
}\]	
\end{proof}

Let $\mathcal{T}$ be an abelian category with enough projective objects and enough injective objects. Let $\mathcal{X}$ and $\mathcal{Y}$ be two full subcategories of $\mathcal{T}$. 
	
\begin{Prop}$($\cite[Proposition 2.3]{IY}$)$\label{pretri-torsion-pair} 
Let $\mathcal{T}=(\mathcal{T},\Sigma,\nabla)$ be a  right triangulated category. Then
\begin{enumerate}	
\item Let  $\mathcal{X}$ be a contravariantly finite {\rm(}resp. covariantly finite{\rm)}  and extension closed subcategories of  $\mathcal{T}$. Then $(\mathcal{X}, \mathcal{X}^{\bot})${\rm(}resp. $({^{\bot}\mathcal{X}},\mathcal{X})${\rm)} forms a torsion theory.
\item  $\mathcal{X}\mapsto\mathcal{Y}:=\mathcal{X}^{\bot}$ gives a one-to-one correspondence between contravariantly finite and extension closed subcategories  $\mathcal{X}$  of $\mathcal{T}$ and covariantly finite and extension closed subcategories  $\mathcal{Y}$  of $\mathcal{T}$. The inverse is given by  $\mathcal{Y}\mapsto\mathcal{X}:={^{\bot}\mathcal{Y}}.$ 
\end{enumerate}	
\end{Prop}
\begin{proof} 

$(1)$ This is a direct generalization of \cite[Proposition 2.3 (1)]{IY}, for the convenience of the reader, we include a proof here.

We only prove that $(\mathcal{X}, \mathcal{X}^{\bot})$  is a torsion pair, the case for $(^{\bot}\mathcal{X}, \mathcal{X})$ can be proved by a dual argument. We know that both $\mathcal{X}$ and $\mathcal{X^{\perp}}$  are  closed under direct summands, and clearly $\mathcal{T}\mathcal(\mathcal{X},\mathcal{X^{\perp}})=0$.
It remains to check the second condition of torsion pair. Let $M\in\mathcal{T}$. Since $\mathcal{X}$ is contravariantly finite in $\mathcal{T}$, we have a
right  triangle
$$\xymatrix@C=0.5cm{
   X \ar[r]^{f} & M\ar[r]^{g} &Y \ar[r]^{} & \Sigma X}$$
with $X\in \mathcal{X}$ and $f$ a right $\mathcal{X}$-approximation of $M$. Without loss of generality, we can assume that $M$ is indecomposable.

By \cite[ Proposition 4.2]{DV}, any morphism $f\colon A\rightarrow B$ in a Krull-Schmidt category (up to isomorphism) has a decomposition $f=(f',0)\colon A'\oplus A''\rightarrow B$ with $f'\colon A'\rightarrow B$ right minimal. It follows that we can choose a minimal version of $f$, which is still denoted by $f$. We finish our proof by showing that $Y \in \mathcal{X^{\perp}}$.

If $M\in \mathcal{X}$, then $Y=0$, if $M\in\mathcal{X^{\perp}}$, then $X=0$ and $Y\cong M\in\mathcal{X^{\perp}}$, so we assume that $M\notin \mathcal{X}\cup\mathcal{X^{\perp}}$.  It is sufficient to show that for any $S\in\mathcal{X}$, $\mathcal{T}(S,Y)=0$. Otherwise, there exists $S_{0}\in\mathcal{S}$ such that  $\mathcal{T}(S_{0},Y)\neq0$. Let $0\neq \rho\in\mathcal{T}(S_{0},Y)$.
We have the following commutative diagram in $\mathcal{T}$:

$$\xymatrix{
   &  & Z \ar[ld]_{\sigma} \ar[d]^{\delta}  &   \\
 X\ar@{=}[d]_{\Id}\ar[r]^{m} & H  \ar[d]_{\varphi} \ar[r]^{n} & S_{0} \ar[d]_{\rho} \ar[rd]^{h\circ\rho} &   \\
  X \ar[r]^{f} & M \ar[r]^{g} & Y \ar[r]^{h} \ar[d]^{} & \Sigma X,  \\
  &   &  \Sigma Z &   \\}$$
where  $\xymatrix@C=0.5cm{Z\ar[r]^{\sigma} & H\ar[r]^{\varphi} & M\ar[r]^{} & \Sigma Z}$
 is the right triangle induced from the octahedral axiom in $\mathcal{T}$. Since $X$ and $S_{0} \in \mathcal{X}$, $H\in\mathcal{X}$. Since $f$ is a right $\mathcal{X}$-approximation of $M$, there exists a morphism $e\colon H\rightarrow X$ such that $\varphi=f\circ e$. Therefore we get a commutative diagram:
\vspace{0.2cm}
$$\xymatrix{
  X  \ar@{=}[d]_{\Id} \ar[r]^{m} & H \ar[ld]_{e} \ar[d]_{\varphi} \ar[r]^{n} & S_{0}\ar[d]^{\rho} \\
  X \ar[r]^{f} & M \ar[r]^{g} & Y. }  $$
\vspace{0.1cm}
\noindent Then $f=\varphi\circ m=f\circ e\circ m$, and since it is minimal, $e\circ m$ is an isomorphism. Therefore we  have that $m$ is a split monomorphism and  $n$ is a split epimorphism. Since $\rho$ is nonzero, $\rho\circ n\neq0$. However we  have that $\rho\circ n=g\circ\varphi=g\circ f\circ e=0$, it is a contradiction. Hence $Y \in \mathcal{X^{\perp}}$.	

$(2)$ The assertion follows from $(1)$.
\end{proof}

\begin{Rem}\label{covariantly-fin-con}
\begin{enumerate}[$(1)$]
\item Dually, our conclusion holds for a left triangulated category. 
\item If $\mathcal{T}$  has  Serre duality and $\mathcal{X}\subseteq\mathcal{T}$ is contravariantly finite or covariantly finite in $\mathcal{T}$, then both  $(\mathcal{X}, \mathcal{X}^{\bot})$ and  $({^{\bot}\mathcal{X}},\mathcal{X})$ are torsion pairs, where Serre duality is an exact auto-equivalence  $\nu:\mathcal{T}\rightarrow\mathcal{T}$ satisfying a natural isomorphisms 
\begin{equation*}\label{serre-dual}
\mathcal{T}(X,Y)\cong \D\mathcal{T}(Y,\nu X)\ for\  all\  X,Y\in\mathcal{T}.
\end{equation*}
Indeed, it is clear that both $(^{\bot}\mathcal{X},(^{\bot}\mathcal{X})^{\bot})$ and  $({^{\bot}(\mathcal{X}^{\bot})},\mathcal{X}^{\bot})$ are torsion pairs.  And, if $\mathcal{T}$  has  Serre duality, for $\mathcal{X}\subseteq\mathcal{T}$, it follows from above isomorphism that   $(^{\bot}\mathcal{X})^{\bot}=\mathcal{X}={^{\bot}(\mathcal{X}^{\bot})}$.
\end{enumerate}
\end{Rem}

From now on, we study the relationship between sms's and wsms's.
Let $R$ be a finite dimensional  $k$-algebra and  $\mathcal{S}$  an orthogonal system in $R$-$\stmod$. Going back  to module category $R$-$\m$ under standard projection $\pi\colon R$-$\m\rightarrow R$-$\stmod$,  we know that  $F(\mathcal{S},R):=\mathcal{F}(\mathcal{S})\cup \add(_RR)$ is the preimage of $\mathcal{F}(\mathcal{S})$ under the projection $\pi$. Note that $F(\mathcal{S},R)$ is closed under extensions and direct summands in $R$-$\m$. Take $M=\bigoplus_{M_{i}\in F(\mathcal{S},R)}M_{i}$ and consider the endomorphism ring $S=\End(_{R}M)$.
We have the following results.

\begin{Them}\label{properties-subset-of-sms}
Let $R$ be a finite dimensional  algebra and $\mathcal{X}$  a subset of an sms $\mathcal{S}$ in $R$-$\stmod$. Let $M=\bigoplus_{M_{i}\in F(\mathcal{X},R)}M_{i}$  and $S=\End(_{R}M)$.  
If $M$ is finendo {\rm(}that is, a finitely generated right $S$-module{\rm)}, then $S$ is  a right coherent ring. 
\end{Them}
\begin{proof}
By Theorem \ref{subset-of-sms}, $F(\mathcal{X})$ is functorially finite in $R$-$\stmod$. Since $\mathcal{P}_{R}=\add(_RR)$ is  covariantly finite in $R$-$\m$,
$F(\mathcal{X},R)$ is covariantly finite in $R$-$\m$ by Lemma \ref{contra-variant}. It follows from  Proposition \ref{direct-prod-definable} that $\varinjlim F(\mathcal{X},R)$ is closed under direct products.
Let $(F_{i}\mid i\in I)$ be  a class of flat  left $S$-modules. It suffices to show that $\prod_{i\in I}F_{i}$
is a flat left $S$-module. 

By  Theorem \ref{direct-lim-of-add}, we may take a class of modules of the form $(M\otimes_{S} F_{i})_{i\in I}$ in $\varinjlim F(\mathcal{X},R)$. Since $\varinjlim F(\mathcal{X},R)$ is closed under direct products, $\prod_{i\in I}M\otimes_{S} F_{i}$ belongs to $\varinjlim F(\mathcal{X},R)$. Since 
$M$ is a finitely generated right $S$-module, by Lemma \ref{coherent-epi-iso},  $M\otimes_{S} \prod_{i\in I}F_{i}\cong\prod_{i\in I}M\otimes_{S} F_{i}$ is in $\varinjlim F(\mathcal{X},R)$. It follows from  Theorem \ref{direct-lim-of-add} that $\prod_{i\in I}F_{i}$ is a flat left $S$-module. Thus $S$ is  a right coherent ring. 
\end{proof}

\begin{Them}\label{fin-dim-sms}
Let $R$ be a finite dimensional   $k$-algebra and $\mathcal{S}$  an orthogonal system in $R$-$\stmod$. Let $M=\bigoplus_{M_{i}\in F(\mathcal{S},R)}M_{i}$  and $S=\End(_{R}M)$.   Then 
\begin{enumerate}[$(1)$]
\item If $S$ is right coherent and $M$ is finendo, then $\mathcal{F}(\mathcal{S})$  is covariantly finite in $R$-$\stmod$. 
\item  If $S$ is right coherent, $M$ is finendo and  $\mathcal{S}$  is a wsms that  contains no  objects of any homogeneous tube,  then $\mathcal{S}$ is an sms. 
\end{enumerate}
\end{Them}
\begin{proof}
(1) We first show that 	$\mathcal{F}(\mathcal{S})$  is covariantly finite in $R$-$\stmod$. Since $\mathcal{P}_{R}$ is functorially finite, by Lemma \ref{contra-variant}, it suffices to show that $F(\mathcal{S},R)$ is covariantly finite in $R$-$\m$. By  Proposition \ref{direct-prod-definable}, $F(\mathcal{S},R)$ is covariantly finite in $R$-$\m$ if and only if $\varinjlim F(\mathcal{S},R)$ is closed under direct products.  Now it suffices to show that $\varinjlim F(\mathcal{S},R)$ is closed under direct products.  

Without loss of generality, by Theorem \ref{direct-lim-of-add}, we may take a class of modules of the form $(M\otimes_{S} Q_{i})_{i\in I}$ in $\varinjlim F(\mathcal{S},R)$, where $(Q_{i}\mid i\in I)$ is a class of flat  left $S$-modules. Since $M_{S}$ is finitely generated, by Lemma \ref{coherent-epi-iso}, we have  the isomorphism 
\[ \varphi_{M,\mathcal{Q}}\colon M\otimes_{S} \prod_{i\in I}Q_{i}\longrightarrow\prod_{i\in I}M\otimes_{S} Q_{i}.\]  
  Since $S$ is left coherent, the product of left flat $S$-modules is left flat, therefore $\prod_{i\in I}Q_{i}$ is still a left flat $S$-module. Hence, by Proposition \ref{direct-lim-of-add}, $\prod_{i\in I}M\otimes_{S} Q_{i}\in \varinjlim F(\mathcal{S},R)$, therefore $\varinjlim F(\mathcal{S},R)$ is closed under direct products. It follows from  Lemma \ref{contra-variant} and Proposition \ref{direct-prod-definable} that $\mathcal{F}(\mathcal{S})$ is covariantly finite in $R$-$\stmod$.
  
(2) We prove that $\mathcal{S}$ is an sms.   If $\mathcal{S}$ is an sms, then by Theorem \ref{simple-module-and-n-tube}, there is no objects in any homogeneous tube contained in $\mathcal{S}$. It is clear that $\mathcal{F}(\mathcal{S})$  is closed under extensions and direct summands in $R$-$\stmod$. By Proposition \ref{pretri-torsion-pair} and Remark \ref{covariantly-fin-con}, $(^{\bot}\mathcal{S}, \mathcal{F}(\mathcal{S}))$ is a torsion pair in  $R$-$\stmod$.  Since $\mathcal{S}$  is a wsms, we have $^{\bot}\mathcal{S}=\{0\}$. Hence $\mathcal{F}(\mathcal{S})=R$-$\stmod$, and  $\mathcal{S}$  is an sms.
\end{proof}

According to  Theorem \ref{coherent-ring-add-kernel-2}  and \ref{coherent-ring-add}, we have the following corollary. 
\begin{Cor}\label{fin-dim-sms-pseuker}
Let $R$ be a finite dimensional  $k$-algebra, $\mathcal{S}$  is a wsms that  contains no  objects of any homogeneous tube in $R$-$\stmod$. Let $M$ and $S$ be as above of Theorem \ref{fin-dim-sms}.  
If $M$ is finendo and $\add(M)$ satisfies one of conditions:   
\begin{enumerate}[$(1)$]
\item  Every  morphism of $\add(M)$ has a pseudokernel,	especially every  morphism of $\add(M)$ has a kernel.
\item Every finitely $M$-cogenerated module has a right  $\add(M)$-approximation, especially, every  finitely $M$-copresented module belongs to $\add(M)$.
\end{enumerate}
Then $\mathcal{S}$ is an sms. 
\end{Cor}

\subsection{Hereditary algebras cases}
Let $A$ be an artin algebra and  $\Gamma(A{\text-}\m)$  the  Auslander-Reiten quiver (AR-quiver for short) of $A$. 
A connected component $\mathcal{P}$ of $\Gamma(A{\text-}\m)$ is called {\bf postprojective} if  $\mathcal{P}$ is acyclic and, for any indecomposable module $M$ in $\mathcal{P}$, there exist $t\geq 0$ and $a\in(Q_{A})_{0}$ such that $M\cong\tau^{-t}P(a)$. 
A connected component $\mathcal{Q}$ of $\Gamma(A{\text-}\m)$ is called {\bf postinjective} if  $\mathcal{Q}$ is acyclic and, for any indecomposable module $N$ in $\mathcal{Q}$, there exist $s\geq 0$ and $b\in(Q_{A})_{0}$ such that $N\cong\tau^{-s}I(b)$. An indecomposable $A$-module is called {\bf postprojective} (resp. {\bf postinjective}) if it belongs to a postprojective component (resp. postinjective component) of $\Gamma(A{\text-}\m)$. 


Recall that an algebra $A$ is said to be {\bf left hereditary} if any left ideal  in $A$ is projective as a left $A$-module.  Note that, for a finite dimensional algebra $A$, $A$ is left hereditary if and only if $A$ is right hereditary, therefore we speak about hereditary without further specification. 

Let $A$ be a hereditary algebra. According to  the known description of AR-quiver $\Gamma(A{\text-}\m)$ of $A$  by Skowro\'{n}ski  \cite{SS}, $\Gamma(A{\text-}\m)$  has the disjoint union form 
\[\Gamma(A{\text-}\m)=\mathcal{P}(A)\cup\mathcal{R}(A)\cup\mathcal{Q}(A).\]
Where 
\begin{enumerate}[$\bullet$]
\item $\mathcal{P}(A)$ is the unique preprojective component containing all the indecomposable projective modules.
	
\medskip
	
\item $\mathcal{Q}(A)$ is the unique preinjective component containing all the indecomposable injective modules.
	
\medskip
	
\item The regular  part $\mathcal{R}(A)$ is a disjoint union of the $\mathbb{P}_{1}(k)$-family  \[{\bf\mathcal{T}}^{A}=\{T_{\lambda}^{A}\}_{\lambda\in\mathbb{P}_{1}(k)}\]
	of stable tubes which are pairwise orthogonal in $A$-$\m$. 
	
Note that at most three tubes in $\mathcal{T}^{A}$ are of rank $\geq2$,  the remaining ones are homogeneous tubes and that  the family $\mathcal{T}^{A}$ separates the postprojective component $\mathcal{P}(A)$ from the preinjective component $\mathcal{Q}(A)$, 
that is, $\Hom_{A}(\mathcal{Q}(A),\mathcal{P}(A))=0$, $\Hom_{A}(\mathcal{Q}(A),\mathcal{R}(A))=0$, and $\Hom_{A}(\mathcal{R}(A),\mathcal{P}(A))=0$.
\end{enumerate}

\bigskip

\begin{center}
	\begin{tikzpicture}
		\draw[]	(0,0) --(0,2);
		\draw[]	(0,2) --(2.5,2);
		\draw[densely dashed]	(2.5,2) --(3.5,2);
		\draw[densely dashed]	(4.5,2) --(5.5,2);
		\draw[]	(5.5,2) --(8,2);
		\draw[]	(8,2) --(8,0);
		\draw[]	(0,0) --(2.5,0);
		\draw[densely dashed]	(2.5,0) --(3.5,0);
		\draw[densely dashed]	(4.5,0) --(5.5,0);
		\draw[]	(5.5,0) --(8,0);
		\draw[] (4,1) circle (1.1);
		\coordinate(1) at (1.5,1);
		\node at (1)[]{$\mathcal{P}(A)$};
		\coordinate(2) at (4,1);
		\node at (2)[]{$\mathcal{R}(A)$};
		\coordinate(3) at (6.5,1);
		\node at (3)[]{$\mathcal{Q}(A)$};
	\end{tikzpicture}
\end{center}

Recall that an algebra $B$ is called {\bf concealed of type $Q$} if there exists a postprojective tilting module $T$ over the path algebra $A=kQ$ such that $B=\End(_{A}T)$, where $Q$ is a finite, connected and acyclic quiver that is not a Dynkin quiver. 
Note that  a hereditary algebra of a Euclidean type is clear a concealed algebra of  Euclidean type by taking left regular $A$-module $_{A}A$ as the postprojective tilting module.

\begin{Lem}$($\cite[Chapter XI, Section 2, Proposition 2.4]{SS}$)$\label{hereditary-extension}
Let $Q$ be an acyclic quiver whose underlying graph is Euclidean and let $A=KQ$. The full subcategory $\add(\mathcal{R}(A))$ of $A$-$\m$ is abelian and closed under extensions. 
\end{Lem}

\begin{Lem}$($\cite[Chapter XII, Section 3, Proposition 3.6]{SS}$)$\label{to-homogenous-morphism}
Let $B$ be a concealed algebra of a Euclidean type $Q$, $\mathcal{T}^{B}=\{\mathcal{T}^{B}_{\lambda}\}_{\lambda\in\mathbb{P}_{1}(k)}$  the family of all stable tubes of $\Gamma(B{\text-}\m)$ and let $M$ be an indecomposable $B$-module. Fix $\lambda\in\mathbb{P}_{1}(k)$ and denote by $r_{\lambda}$ the  rank of the tube $\mathcal{T}^{B}_{\lambda}$. 
\begin{enumerate}[$(1)$]
\item If $M$ is postprojective, then $\Hom_{B}(M,X)\neq0,$ for any indecomposable $B$-module $X$ in the tube $\mathcal{T}^{B}_{\lambda}$ such that $r\ell(X)\geqslant r_{\lambda}$. 
		
\item If $M$ is postinjective, then $\Hom_{B}(X,M)\neq0,$ for any indecomposable $B$-module $X$ in the tube $\mathcal{T}^{B}_{\lambda}$ such that $r\ell(X)\geqslant r_{\lambda}$. 
	\end{enumerate}
\end{Lem}	

\begin{Prop}$($\cite[Chapter 8, Corollary 8.43]{GT}$)$\label{covariantly-finite-1} 
Let $A$ be a left hereditary and right noetherian ring. Assume that $\mathcal{C}\subseteq A$-$\m$ is closed under extensions and direct summands and $A\in\mathcal{C}$.  Then $\mathcal{C}$ is covariantly finite in $A$-$\m$.
\end{Prop}

\begin{Cor}\label{sms-obj}
Let $Q$ be an acyclic quiver whose underlying graph is Euclidean and let $A=KQ$. There is at least one object of the postprojective component or the postinjective component contained in an sms over $A$-$\stmod$.
\end{Cor}
\begin{proof}
It is a direct consequence of Lemma \ref{hereditary-extension}.
\end{proof}

\begin{Prop}\label{hereditary-homogene-sms}
Let $A$ be a finite dimensional hereditary algebra of Euclidean type. Then there is no objects in any homogeneous tube contained in an sms.
\end{Prop}
\begin{proof}
We know that a  hereditary algebra of Euclidean type is always a concealed algebra of  Euclidean type, thus it is sufficient to show that our conclusion holds for  a concealed algebra of  Euclidean type.  According to  the known description of AR-quiver $\Gamma(A{\text-}\m)$ of $A$  by Skowro\'{n}ski \cite{SS}, $\Gamma(A{\text-}\m)$ consists of the unique preprojective component  $\mathcal{P}(A)$ containing all the indecomposable projective modules, the unique preinjective component $\mathcal{Q}(A)$ containing all the indecomposable injective modules and a family of regular components $\mathcal{R}(A)$ containing a $\mathbb{P}_{1}(k)$-family  ${\bf\mathcal{T}}^{A}=\{T_{\lambda}^{A}\}_{\lambda\in\mathbb{P}_{1}(k)}$ of stable tubes which are pairwise orthogonal in $A$-$\m$. At most three tubes in $\bf{\mathcal{T}}$ are of rank $\geq2$,  the remaining ones are homogeneous tubes.
	
By Lemma \ref{to-homogenous-morphism}, it is clear that, for an indecomposable postprojective $A$-module $M$, $\Hom_{A}(M,X)$ is always not zero space  for any module $X$ in a homogeneous tube. Dually, for an indecomposable postinjective $A$-module $M$, $\Hom_{A}(X,M)$ is always not zero for any module $X$ in a homogeneous tube. By Corollary \ref{sms-obj},  there exists at least one object coming from the postprojective component or the postinjective component in an sms. By the orthogonality of sms's,  there is no objects in any homogeneous tube contained in an sms.
	
\end{proof}

\begin{Them}\label{weak-sms}
Let $A$ be a finite dimensional hereditary $k$-algebra and $\mathcal{S}$  an orthogonal system in $A$-$\stmod$, then $\mathcal{F}(\mathcal{S})$ is covariantly finite in $A$-$\stmod$. If moreover $\mathcal{S}$  is a weakly simple-minded system with finite cardinality in $A$-$\stmod$, then $\mathcal{S}$ is a simple-minded system in $A$-$\stmod$.
\end{Them}
\begin{proof}
We know that $F(\mathcal{S},A)$ is closed under extensions and direct summands.  By Proposition \ref{covariantly-finite-1}, $F(\mathcal{S},A)$ is covariantly finite in $A$-$\m$. By Lemma \ref{contra-variant}, $\mathcal{F}(\mathcal{S})=F(\mathcal{S},A)/\mathcal{P}_{A}$  is covariantly finite in $A$-$\stmod$. 
Moreover, if $\mathcal{S}$  is a weakly simple-minded system with finite cardinality in $A$-$\stmod$, then by   Proposition  \ref{hereditary-homogene-sms}, $\mathcal{S}$  contains no objects of any homogeneous tube. 
It follows from Proposition \ref{pretri-torsion-pair} that  $(^{\bot}\mathcal{S},\mathcal{F}(\mathcal{S}))$ is a torsion pair in $A$-$\stmod$. Since $^{\bot}\mathcal{S}=\{0\}$, $\mathcal{F}(\mathcal{S})=A$-$\stmod$. Thus $\mathcal{S}$ is a simple-minded system in $A$-$\stmod$.
\end{proof}

\subsection{pure-semisimple ring case}
Recall that a ring $A$ is said to be left (resp. right) {\bf pure-semisimple} if every left (resp. right) $A$-module is a direct sum of finitely presented modules, or equivalent, if every left (resp. right) $A$-module is pure-injective. Note that a $A$-module $M$ is called pure-injective if it has the injective property relative to the class of pure-exact sequences. 
It is known that a ring is of finite representation type if and only if it is left and right pure-semisimple. However, it has been a long-standing open problem, known as the {\bf pure semisimplicity conjecture}, whether left  pure-semisimple rings are also right pure semisimple (see  for example \cite{ZH2}). 
Finendo modules play a key role when we consider the relationship between sms's and coherent rings (see Theorem \ref{properties-subset-of-sms} and Theorem \ref{fin-dim-sms}). The following conclusion indicates a close relationship between finendo modules and pure semisimple rings.

\begin{Prop}{\rm(}\cite[Theorem 2.11 and Corollary 2.12]{DG}{\rm)}\label{finendo-pure-semisimple}
Let $A$ be a  ring. If every right $A$-module is finendo, then $A$ is left pure semisimple and every finitely generated indecomposable left $A$-module is cofinendo.
\end{Prop}

Under hereditary assumption, there is a further relationship as follows.

\begin{Prop}{\rm(}\cite[Theorem 3.5]{DG}{\rm)}\label{hereditary-pure-semisimple}
	Let $A$ be a hereditary ring.  Then the following conditions are equivalent.
\end{Prop}
\begin{enumerate}
\item Every right $A$-module is finendo.
\item $A$ is left pure semisimple, and every finitely generated indecomposable left $A$-module is finendo.
\item $A$ is of finite representation type.
\end{enumerate}

Under pure-semisimple assumption, we have the following well-behaved results. 

\begin{Prop}{\rm(}\cite[Proposition 6.17]{AH}{\rm)}\label{pure-semisimple}
Let $A$ be a left pure-semisimple ring.  Then every class of modules $\mathcal{M}$ is covariantly finite in $A$-$\m$.
\end{Prop}
\begin{proof}
It is a direct consequence of \cite[Theorem 9]{ZH} and \cite[Corollary 3.3]{AH}.
\end{proof}

\begin{Cor}\label{pure-semisimple-sms}
Let $A$ be a left pure-semisimple ring and $\mathcal{S}$  an orthogonal system on  $A$-$\stmod$. If $\mathcal{S}$ is a wsms, then $\mathcal{S}$ is an sms.
\end{Cor}
\begin{proof}
	It is a direct consequence of Theorem \ref{fin-dim-sms} and Theorem \ref{subset-of-sms}.
\end{proof}

\section{Declarations}

\subsection*{ Ethical Approval}

This declaration is “not applicable”.

\subsection*{ Funding} 

The research work is supported by NSFC (No.12031014 and No.12301044).

\subsection*{Availability of data and meterials}
Data availability is not applicable to this article as no new data were created or analyzed in this study.

\end{document}